\documentclass[12pt,a4paper]{article}
\usepackage{tikz}
\usetikzlibrary{shapes.geometric, arrows}

\usepackage{amsmath,amssymb}
\usepackage{mathtools}
\usepackage{etex}
\usepackage{bm}
\usepackage{color}
\newcommand{\Rx}{\mathbb{R}}
\newcommand{\teta}{\tilde{\eta}}
\newcommand{\tim}{\tilde{m}}

\newcommand{\half}{\small{\frac{1}{2}}}
\newcommand{\Px}{\mathbb{P}}

\renewcommand{\span}{{\rm span}}

\newcommand{\uq}{\hat{q}}

\newcommand{\uni}{\underline{1}}
\newcommand{\qed}{\hfill$\Box$}

\title{{\Large \bf Optimal 
approximations of the Fokker Planck Kolmogorov equation: projection, maximum likelihood eigenfunctions and Galerkin methods}}
\author{Damiano Brigo \\ Dept. of Mathematics \\ Imperial College London \\ 180 Queen's Gate \\  {\small \tt{damiano.brigo@imperial.ac.uk}}\and   Giovanni Pistone \\de Castro Statistics \\ Collegio Carlo Alberto \\ Via Real Collegio 30 \\ 10024 Moncalieri, IT }

\newtheorem{theorem}{Theorem}[section]
\newtheorem{proposition}[theorem]{Proposition}

\newtheorem{remark}[theorem]{Remark}

\begin{document}
%
%\PARstart{X}{YYY} ZZZ

\newpage

\maketitle

%
%
%

%
%\PARstart{X}{YYY} ZZZ

%\thispagestyle{empty}
%
%
%

\begin{abstract}
%We apply the $L^2$ based Hellinger-Fisher-Rao vector-field projection introduced in Brigo, Hanzon and LeGland  (1998, 1999) to finding locally and globally optimal finite dimensional approximations of the Fokker--Planck or forward Kolmogorov equation on exponential families. We show that this vector field projection is equivalent to a local assumed density ``moment matching'' approximation based on expectation parameters, in line with a previous more general result for the stochastic PDE of the filtering problem. We derive an algebraic relation that allows to recover canonical parameters from expectation parameters for polynomial exponent families. The main result of this paper is showing that, for general exponential families, if the sufficient statistics are chosen among the diffusion eigenfunctions then the finite dimensional projection or the equivalent assumed density approximation provide the exact maximum likelihood density and we have global optimality in relative entropy. Applications include the prediction step for the filtering problem with discrete time observations. Our result is based on the differential geometric approach to statistics and systems theory.   
%
%
%
We study optimal finite dimensional approximations of the generally infinite-dimensional Fokker-Planck-Kolmogorov (FPK) equation, finding the curve in a given finite-dimensional family that best approximates the exact solution evolution. For a first local approximation we assign a manifold structure to the  family and a metric. We then  project the vector field of the partial differential equation (PDE) onto the tangent space of the chosen family, thus obtaining an ordinary differential equation for the family parameter. A second global approximation will be based on projecting directly the exact solution from its infinite dimensional space to the chosen family using the nonlinear metric projection. This will result in matching expectations with respect to the exact and approximating densities for particular functions associated with the chosen family, but this will require knowledge of the exact solution of FPK.  A first way around this is a localized version of the metric projection based on the assumed density approximation. While the localization will remove global optimality, we will show that the somewhat arbitrary assumed density approximation is equivalent to the mathematically rigorous vector field projection. More interestingly we study the case where the approximating family is defined based on a number of eigenfunctions of the exact equation. In this case we show that the local vector field projection provides also the globally optimal approximation in metric projection, and for some families this coincides with a Galerkin method. 
We study exponential and mixture families, and the metrics for the vector field projection are respectively the Hellinger-Fisher-Rao and the direct $L^2$ distances. For the metric projection we use respectively relative entropy and $L^2$ direct distance. In the eigenfunctions case we derive the exact maximum likelihood density for FPK. Our results are based on the differential geometric approach to statistics and systems theory, and applications include filtering. 

\end{abstract}

\medskip


\medskip

% REQUIRED

\newpage

\tableofcontents

\vfill

%\vspace{6cm}

\subsection*{Acknowledgments}
{\small The authors are grateful to John Armstrong for many stimulating and interesting discussions and for geometric intuition that helped improve the paper. In particular, the metric projection interpretation of the global approximation is based on an initial suggestion by John in related work for optimal approximation of S(P)DEs on submanifolds (\cite{armstrongBrigoproj} and \cite{armstrongbrigoicms}, building on \cite{armstrongBrigoJets}).}

\pagestyle{myheadings}
\markboth{}
{D. Brigo, G. Pistone. Optimal approximation of the Fokker-Planck-Kolmogorov eq. }

\newpage

\section{Introduction}

Problems in systems theory, especially filtering and control, may have solutions that are expressed as evolutions of probability distributions. Approximating the evolution of a probability density function with an evolution in a parametric family is an important approximation problem. If one is working in a setting where a probability density $p_t(x)$ at the point $x$ for every time $t \ge 0$ evolves in time as a curve $t \mapsto p_t$ on an infinite dimensional space, it may be important to have a sound methodology to find the best approximation of the curve $t \mapsto p_t$ in a finite dimensional parametrized family of densities, say $\{ p(\cdot,\theta), \theta \in \Theta \subset \Rx^n\}$. This way the approximated density $p(x,\theta_t)$ will be a natural approximation of the full density $p_t(x)$ and it will be much more manageable given the dimensionality reduction from infinite to a finite $n$. If one has an evolution equation for $p_t$, say a partial differential equation (PDE) or a Stochastic PDE, as happens for example in the filtering problem with the Fokker-Planck-Kolmogorov (discrete time observations) or Kushner-Stratonovich / Zakai (continuous time observations) equations, then it becomes very important to approximate such (S)PDE with an equation for $\theta_t$, which will be a (stochastic) differential equation in the finite dimensional space $\Rx^n$. Indeed, any implementation of an equation on a machine needs to be finite dimensional, so that finding optimal finite dimensional approximations $t \mapsto \theta_t$ of the original infinite-dimensional solution $t \mapsto p_t$ is of great practical importance. 
For the filtering problem this has been addressed in the work initially sketched in Hanzon (1987) \cite{hanzon87} and fully developed in  Brigo et al (1998,1999)\cite{brigo98, brigo99}, where a general method for deriving finite dimensional approximations of the optimal filter stochastic PDE in a chosen exponential family has been given.  
In this paper we focus on the simpler Fokker-Planck-Kolmogorov (FPK) equation. We still have a strong link with the filtering problem, since under discrete time observations the FPK equation constitutes the prediction step between observations and is still a fundamental part of the filtering algorithm, see for example \cite{jazwinski70a}.
Applications of finite dimensional approximations of the FPK equation are by no means limited to the filtering problem. Examples of possible applications include, beside signal processing, stochastic- (local-) volatility modeling in quantitative finance, the anisotropic heat equation in physics, and quantum theory evolution equations, see the related discussion and references in \cite{brigopistone}. 
In order to make sense of the notion of ``best'' or ``optimal'' approximation, we need to have a measure of how close an approximation will be to the exact solution, and aim to find the closest. Mathematically, this takes the form of a metric or possibly a divergence in the space of all possible probability densities $p$. This way we will have a notion of distance of $p(\cdot,\theta_t)$ from $p_t$, a distance we may wish to minimize in some sense to find the best possible approximation. ``Optimality'' can be imposed either in a local sense, or more strongly in a global sense. Local optimality means that whenever we evolve away {\emph{locally}} from the finite dimensional family $\{ p(\cdot,\theta), \theta \in \Theta \}$ following the ``vector field'' of the original infinite dimensional (S)PDE, resulting in a ``$d {p_t}|_{p_t=p(\cdot, \theta_t)}$'' that points out of the (tangent space of the) family,  we will find the vector $d p(\cdot, \theta_t)$ staying in the (tangent space of the) family that is closest to $d p_t$, and follow $d p(\cdot, \theta_t)$ rather than the full infinite-dimensional evolution. This is only local optimality since we approximate a vector departing from our family with a vector staying in our family, but we don't approximate directly the solution, which leaves the family immediately. The main technique used for this local optimality will be, not surprisingly, the linear projection (on the tangent space), and the metric we will use in the space of densities will be mostly the $L^2$ metric on square root of densities $\sqrt{p}$, the Hellinger distance $d_H(p_1,p_2)^2 = \int (\sqrt{p_1(x)} - \sqrt{p_2(x)})^2 dx$, or alternatively the $L^2$ distance taken directly on densities $p$ themselves $d_D(p_1,p_2)^2 = \int (p_1(x) - {p_2(x)})^2 dx$, assuming these are square integrable. We will refer shortly to the projection on the tangent space as to ``tangent'' or ``linear projection". 

Global optimality will clearly be stronger than local optimality. For global optimality we wish to find for every time $t$ the density $p(\cdot,\theta_t)$ in the finite dimensional family that is closest to the true solution $p_t$ in a metric or divergence defined in the infinite dimensional space where $p_t$ evolves. When the approximation $p(\cdot,\theta_t)$ evolves in an exponential family $p(x,\theta) = \exp\left(\theta_1 c_1(x) + \theta_2 c_2(x) + \ldots +\theta_n c_n(x) - \psi(\theta) \right)$
we will find the globally optimal approximation by minimizing the relative entropy or Kullback Leibler information with respect to $\theta$, $K(p_t,p(\cdot,\theta) )= \int \ln(p_t(x)/p(x,\theta)) p_t(x) dx$, rather than the Hellinger distance, because this turns out to be much easier when the family $p(\cdot, \theta)$ is an exponential family. This is a sort of metric projection of $p_t$ onto the chosen family $p(\cdot,\theta)$ in relative entropy (not exactly a metric projection because relative entropy is a divergence and not a metric). We will refer to this projection as to the ``metric projection''. This projection is nonlinear in general, although its linearization leads precisely to the tangent projection. When applied to our problem, this metric projection approach will give us a ``moment matching" characterization of the maximum likelihood exponential density, or of the globally optimal approximation. The two are the same because the maximum likelihood density results indeed from minimization of relative entropy. We will show that the globally optimal or maximum likelihood exponential family density approximating $p_t$ will be the one sharing the expectations of the chosen family sufficient statistics with $p_t$. 

These expectations provide another parameterization of the exponential family, alternative to the canonical parameters $\theta$. The expectation parameters $\eta$ are readily computed given the $\theta$ via their definition, which is $\eta_i(\theta) = \int c_i(x) p(x,\theta) dx = \partial_{\theta_i} \psi$. In this paper we introduce an algebraic relation and algorithm to invert this transformation and obtain the $\theta$ given the $\eta$ in the scalar state space case with monomials sufficient statistics $p(x,\theta) = \exp\left(\theta_1 x + \theta_2 x^2 + \ldots \theta_n x^n - \psi(\theta) \right)$, summarizing earlier results in Brigo et al (1996, 1998) \cite{BrigoPhD,brigoime}. 

A localized version of the globally optimal projection, based on the assumed density approximation and on the expectation parameters $\eta$, will turn out to be identical to the local ``vector field'' $d_H$ optimality above, and this is explained by the fact that relative entropy and the Hellinger distance coincide at the lowest order of approximation. 

For the assumed density approximation in filtering we refer to \cite{kushner}, \cite{maybeck} for the approximation with Gaussian densities, and \cite{brigo99} for the more general approximation with exponential families. The equivalence between the Hellinger projection and the assumed density approximation we present here for the FPK  equation and exponential families is a special case of the equivalence for the filtering stochastic partial differential equation presented in \cite{brigo99}. In the earlier reference \cite{hut} this equivalence had been established for the Gaussian family. 

Summarizing the above results in a nutshell, given an exponential family, 

\medskip

\begin{center}
globally optimal approximation via metric projection in relative entropy \\=\\ maximum likelihood estimation \\=\\ moment matching for sufficient statistics\\ \hspace{1cm} \\
and\\ \hspace{1cm} \\
localized version of globally optimal approx via assumed density approximation 
%local (assumed density) moment matching for sufficient statistics 
\\=\\ locally optimal approximation \\=\\ tangent (vector field) projection in Hellinger/Fisher-Rao metric. 
\end{center}

\medskip

Still, even with these equivalence results proven, finding the fully optimal global approximation is hard. Solving the optimal problem is in general impossible without knowing the true density $p_t$. However, we will state a theorem where for special exponential families $p(\cdot,\theta)$ the local and globally optimal approximations coincide and can be computed without knowing the true infinite-dimensional solution $p_t$. This special case is the main result of the paper and will be based on an exponential family built on the eigenfunctions of the operator associated with the original infinite dimensional equation for the full curve $t \mapsto p_t$. In a nutshell:

\medskip

\begin{center}
sufficient statistics chosen among eigenfunctions of FPK equation operator 
\\ $\Downarrow$\\
 the locally optimal projection is also globally optimal and provides maximum likelihood estimator.
\end{center}

\medskip

A second new series of results of the paper is that we will also show analogous statements to hold for simple mixture families: in that case, to find a globally optimal approximation, our metric projection will be based on minimizing the $L^2$ direct distance $d_D$ between $p_t$ and $p(\cdot,\theta)$. We will show that this also corresponds to an adjusted moment matching condition, and that a localized version of it based on the assumed density approximation is also equivalent to the tangent ``vector-field'' projection of the locally optimal approximation in $d_D$ direct metric. We also recover equivalence of the two with Galerkin type methods. Equivalent to the Galerkin method was shown more generally for the filtering SPDE in \cite{armstrongbrigomcss}.

This can be briefly summarized as follows: given a simple mixture family, 

\medskip

\begin{center}
globally optimal approximation via metric projection in $L^2$ direct distance $d_D$\\=\\  moment matching for mixture components \\ \hspace{1cm} \\
and\\ \hspace{1cm} \\
localized version of globally optimal approx via assumed density approximation 
%local (assumed density) moment matching for sufficient statistics 
\\=\\ locally optimal approximation \\=\\ tangent (vector field) projection in $L^2$ direct metric \\=\\ Galerkin method. 
\end{center}

\medskip

Finally, also for the mixture/$d_D$ case we will show that in the special case where the mixture components are chosen among the FPK operator eigenfunctions, we obtain that the locally optimal tangent projection coincides with the globally optimal metric projection and that local and global optimality coincide.  Again in a nutshell:

\medskip

\begin{center}
mixture components chosen among eigenfunctions of FPK equation operator 
\\ $\Downarrow$\\
 the locally optimal projection is also globally optimal and provides the metric projection.
\end{center}

\medskip

%
%We can summarize the different choices in Table \ref{table:choicesdp}.
%
%\begin{table}
%\begin{tabular}{|c|cccc|}
%\hline 
% Vector    &     &   &   &\\
%Field  &  family: &  Exponential & & Simple mixture \\ 
%   Projection &  &  & & \\ \hline 
%  Metric:      &   &  &  &\\
%  $d_H$     & &  $\sim$Assumed Density &  & \\
%  (Hellinger)              & & Approximation & & \\ 
%                               & & (MLE eigenfunctions) & & \\
%                 & & & & \\
%  $d_D$     & &          & & $\sim$Galerkin Method \\
%  (Direct) $L^2$ & & & & \\
%  \hline
%\end{tabular}
%\caption{Equivalence of projection}\label{table:choicesdp}
%\end{table}   
%
%
%

%
% Add d_D projection on mixtures = Galerkin for Fokker Planck? 
%

\begin{figure}
\begin{tikzpicture}[node distance=2cm]
\tikzstyle{startstop} = [rectangle, rounded corners, minimum width=3cm, minimum height=1cm,text centered, draw=black]
\tikzstyle{arrow} = [thick,->,>=stealth]
\tikzstyle{arrow2} = [thick,<->,>=stealth]
\node (start)[startstop] {Local optimality - tangent projection};
\node (hellinger) [startstop, below of=start,xshift=-2cm] {Hellinger, $d_H$};
\node (direct) [startstop, right of=hellinger,xshift=2cm] {direct, $d_D$};
\node (ef) [left of=hellinger,xshift=-2.3cm] {Exp. Families (EF$(c)$)} ;
\node (sm) [right of=direct,xshift=2.3cm] {Simple Mix. (SM$(\uq)$)} ;
\node (efadf) [startstop, below of=hellinger] {EF-ADA};
\node (smadf) [startstop, below of=direct] {SM-ADA};
\node (expmatchef) [startstop, below of=efadf] {$\mathbb{E}[c]$ matching};
\node (relentropy) [left of=expmatchef,xshift=-2cm]{Relative entropy};
\node (smmatchsm) [startstop, below of=smadf] {$\mathbb{E}[q]$ matching};
\node (directmetric) [right of=smmatchsm,xshift=2cm]{Direct $d_D$};
\node (global) [startstop, below of=smmatchsm,xshift=-2cm]{Global optimality - metric projection};
\draw [arrow] (start) -- (hellinger);
\draw [arrow] (start) -- (direct);
\draw [arrow] (ef) -- (hellinger);
\draw [arrow] (sm) -- (direct);
\draw [arrow2] (hellinger) -- (efadf);
\draw [arrow2] (direct) -- (smadf);
\draw [arrow] (expmatchef) -- node {localization} (efadf);
\draw [arrow] (smmatchsm) -- node {localization} (smadf);
\draw [arrow] (global) -- (expmatchef);
\draw [arrow] (global) -- (smmatchsm);
\draw [arrow] (relentropy) -- (expmatchef);
\draw [arrow] (global) -- (smmatchsm);
\draw [arrow] (directmetric) -- (smmatchsm);
\end{tikzpicture}
\caption{Summary of the different approximations we will introduce in this paper, their optimality and their relationship. If the functions $c$ and $\uq$ are chosen among the eigenfunctions of the given FPK equation, then the locally optimal approximations are also globally optimal.}\label{fig:flow}
\end{figure}
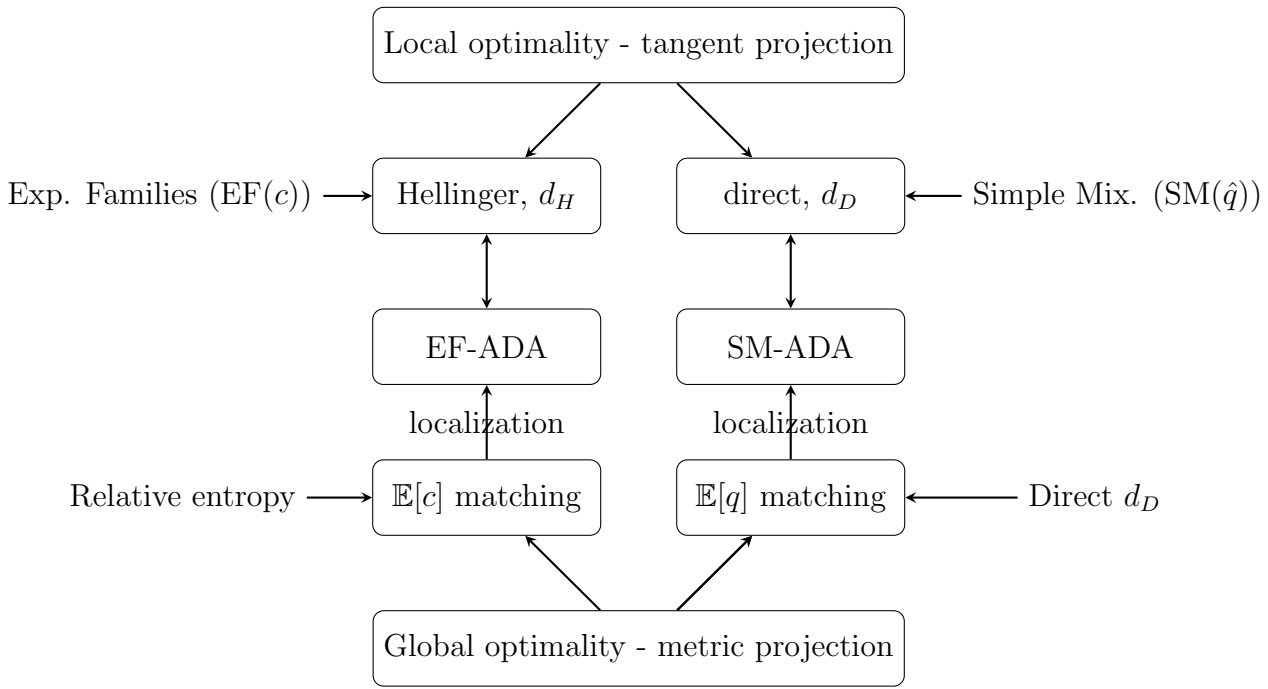

For further references and a detailed literature review see the proceedings paper \cite{brigopistone}, where the maximum likelihood eigenfunctions result for the exponential families case is presented under the different statistical manifolds geometry of Pistone and Sempi \cite{pistonesempi}, based on Orlicz spaces and charts, rather than on the minimal  $L^2$ structure we use here. For general approaches that combine the $L^2$ geometry used here and the Orlicz-based geometry with applications to filtering see for example \cite{Newton1,Newton2,Newton3}.  Notice however that our earlier proceedings paper \cite{brigopistone} does not provide conditions for the existence of the solution of the original equation in the given function space, contrary to our existence result for the $L^2$ structure here, but uses the $L^2$ case itself to proceed, so that the present paper presents the only fully rigorous and consistent analysis on the eigenfunctions maximum-likelihood theorem, see the $L^2$ existence discussion in Section~\ref{sec:vectorfieldprojhellinger} in particular. This is based again on the Hellinger distance. Moreover, this paper deals with the direct distance case as well, which was absent in the proceedings paper. More generally, measures in the space of probability distributions have been used effectively also in sampling high dimensional and strongly correlated systems, see for example the work on Hamiltonian and Langevin Monte Carlo sampling in \cite{girolami}.
We conclude this introduction by mentioning that part of this paper appeared previously as a preprint  in \cite{brigopistone2}. Finally, Figure \ref{fig:flow} summarizes the relationship we'll show in the paper between different approximation methods. 

\bigskip

\section{Statistical manifolds} \label{StatMan}
For a full summary see \cite{brigo99}. We consider parametric families of probability densities, $\{p(\cdot,\theta), \ \theta \in \Theta\}$ with $\Theta$ convex open set in $\Rx^n$. The set of square roots of such densities is a subset of $L^2$ that we may view as a finite dimensional manifold. In general the $L^2$ distance between square roots of densities leads to the Hellinger distance. A curve in such a manifold is given by $t \mapsto \sqrt{p(\cdot,\theta(t))}$. Differentiating with respect to $t$ we obtain that all tangent vectors at $\theta$ are in the space 
$\span\left\{ \frac{\partial{\sqrt{p(\cdot,\theta)}}}{\partial \theta_i}, i=1,\ldots,n\right\}$. 
We can use the $L^2$ inner product $\langle \cdot,\cdot \rangle$ to introduce an inner product on the tangent space and a metric. Recall that for $f,h \in L^2$ we have
$\langle f,h \rangle = \int f(x) h(x) dx$ and recall the related $L^2$ norm
$\| f \|_2 = \left(\int f(x)^2 dx\right)^{1/2}$. 
Define $ g_{i,j}(\theta)/4 = \left\langle \frac{\partial{\sqrt{p(\cdot,\theta)}}}{\partial \theta_i}, \frac{\partial{\sqrt{p(\cdot,\theta)}}}{\partial \theta_j}\right \rangle$. This is, up to the factor $4$, the familiar Fisher-Rao information matrix. 
If we have a $L^2$ vector $v$, we can project it via the $d_H$ tangent or linear orthogonal projection
\begin{equation} \label{proL2}
\Pi^g_\theta[v] = \displaystyle \sum_{i=1}^n [ \sum_{j=1}^n
   4 g^{ij}(\theta)\; \langle v,\frac{1}{2 \sqrt{p(\cdot,\theta)}}\,
   \frac{\partial p(\cdot,\theta)}{\partial \theta_j} \rangle ]\;
   \frac{1}{2 \sqrt{p(\cdot,\theta)}}\,
   \frac{\partial p(\cdot,\theta)}{\partial \theta_i}\ 
\end{equation}
where upper indices denote the inverse matrix. 
A different possibility is a geometry that does not use the square root. This would be done most generally using a duality argument involving $L^1$ and $L^\infty$, but here we will assume that all densities are square integrable, so that all $p$'s we deal with are in $L^2$. Then we can mimic the above structure but without square roots. 
The curve is $t \mapsto p(\cdot,\theta(t))$, the tangent space is 
$\span\left\{ \frac{\partial{{p(\cdot,\theta)}}}{\partial \theta_i}, i=1,\ldots,n\right\}$. 
Define $\gamma_{i,j}(\theta) = \left\langle \frac{\partial{{p(\cdot,\theta)}}}{\partial \theta_i}, \frac{{\partial {p(\cdot,\theta)}}}{\partial \theta_j}\right \rangle$. This leads to what we called $d_D$ tangent or linear projection in \cite{armstrongbrigomcss},
$\Pi^\gamma_\theta[v] = \displaystyle \sum_{i=1}^n [ \sum_{j=1}^n
   \gamma^{ij}(\theta)\; \langle v,
   \frac{\partial p(\cdot,\theta)}{\partial \theta_j} \rangle ]\;
   \frac{\partial p(\cdot,\theta)}{\partial \theta_i}$.
There is another way of measuring how close two densities are.
Consider the Kullback--Leibler information or relative entropy between two densities $p$
and $q$: $K(p,q) := \int \log \frac{p(x)}{q(x)}\; p(x)\, d x$.
This is not a metric, since it is not symmetric and it does
not satisfy the triangular inequality. 
It is a classic result that the Fisher metric and the Kullback--Leibler
information coincide infinitesimally.
Indeed, by Taylor expansion it is easy to show that 
\begin{eqnarray}\label{eq:klig}
 K(p(\cdot,\theta),p(\cdot,\theta+d\theta))
   = \sum_{i,j=1}^n g_{ij}(\theta)\, d\theta_i\, d\theta_j
   + O(\vert d\theta \vert^3)\ .
\end{eqnarray}

\section{Locally optimal approximations to FPK}
In this section we study locally optimal approximations to the FPK equation solution on a family $\{p(\cdot,\theta), \ \ \theta \in \Theta\}$. Our main tool will be projection on the tangent space of the family.

\subsection{General tangent linear (vector-field) projection of FPK}
We first summarize the key results for tangent projection in Hellinger distance $g$ on exponential families and for tangent projection in direct metric $\gamma$ on mixture families. For a full account see \cite{brigo99} for the Hellinger case and \cite{armstrongbrigomcss} for the direct metric case. 
Consider a stochastic differential equation on a probability space $(\Omega, {\cal F}, ({\cal F}_t)_t, \Px)$ taking values in $\Rx^N$,  
\[ dX_t = f(X_t,t) dt + \sigma(X_t,t) dW_t , \ \  X_0, \ \ a(x,t) = \sigma(x,t)\sigma'(x,t)  \]
where the prime index denotes transposition. Let us start from the FPK equation for the probability density $p_t(x)=p(x,t)$ of the solution $X_t$ of our SDE. Examples of possible applications of approximating the FPK equation and its stochastic PDE extensions are given in \cite{brigopistone}, and include signal processing, stochastic- (local-) volatility modeling in quantitative finance, the anisotropic heat equation in physics, and quantum theory evolution equations. 

For our SDE the FPK equation reads
\begin{displaymath}
   \frac{\partial p(x,t)}{\partial t} = {\cal L}_t^\ast\, p(x,t), \ \  {\cal L}_t^\ast p = -\sum_{i=1}^N
   \frac{\partial}{\partial x_i}\, [f_i(\cdot,t)\, p]
   + \half \sum_{i,j=1}^N
   \frac{\partial^2}{\partial x_i \partial x_j}\,
   [a_{ij}(\cdot,t)\, p].
\end{displaymath}
%
% \begin{equation} \label{Kushner}
%    dp_t = {\cal L}_t^\ast\, p_t\,dt
%    + \sum_{k=1}^d p_t\, [\gamma_t^k-E_{p_t}\{\gamma_t^k\}]\,
%    [dY_t^k-E_{p_t}\{\gamma_t^k\}\,dt]
% \end{equation}
% where $E_{p_t}\{\cdot\}$ denotes the expectation w.r.t.\
% the probability density $p_t$, i.e.\ the conditional expectation
% given the observations up to time $t$,
% and where ${\cal L}_t^\ast$ is the forward diffusion operator
% given by
% \begin{displaymath}
%    {\cal L}_t^\ast \phi = -\sum_{i=1}^n
%    \frac{\partial}{\partial x_i}\, [f_t^i\, \phi]
%    + \half \sum_{i,j=1}^n
%    \frac{\partial^2}{\partial x_i \partial x_j}\,
%    [a_t^{ij}\, \phi]\ ,
% \end{displaymath}
% for any test function $\phi$ defined on ${\Rx}^n$,
% i.e.\ ${\cal L}_t^\ast$ is the adjoint operator of
% \begin{displaymath}
%    {\cal L}_t = \sum_{i=1}^n f_t^i\,
%    \frac{\partial}{\partial x_i} + \half \sum_{i,j=1}^n
%    a_t^{ij}\, \frac{\partial^2}{\partial x_i \partial x_j}\ .
% \end{displaymath}
%

\begin{figure}[htp]
\begin{centering}
\includegraphics[scale=0.5]{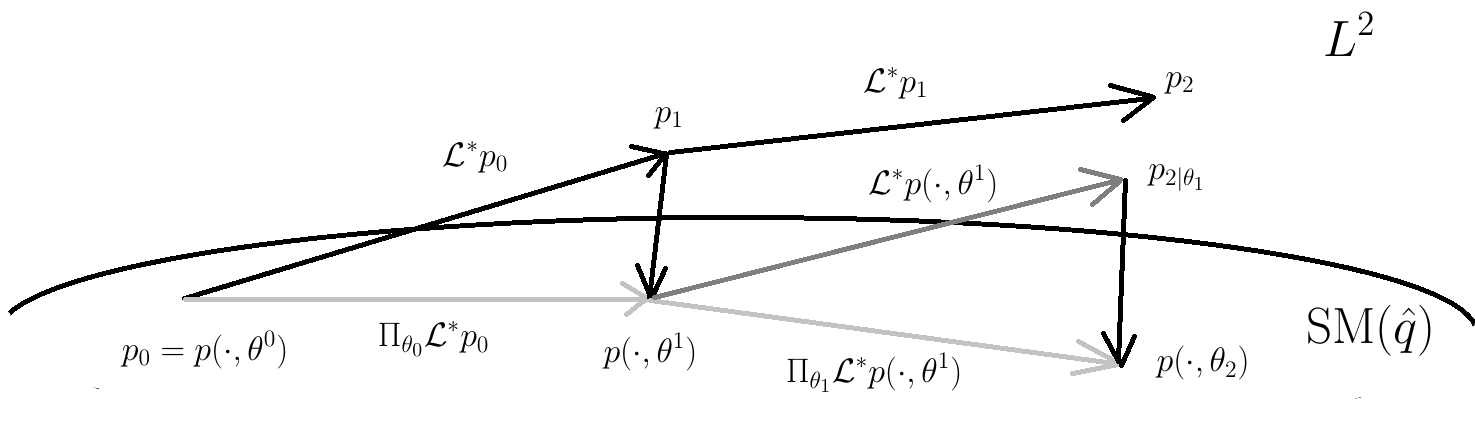}
\caption{Vector field (tangent linear) projection of the FPK Equation on the manifold SM$(\uq)$ in direct metric. We start from an initial condition $p_0$ in the manifold. $p_0,p_1,p_2$ describes the exact solution in $L^2$, whereas $p_0,p(\cdot,\theta^1), p(\cdot,\theta^2)$ describes the evolution resulting from the projected vector field. }
\label{fig:projections}
\end{centering}
\end{figure}

Even if the state space of the underlying diffusion $X$ is finite dimensional, in many cases the FPK equation for the density of $X$ is infinite dimensional (see for example \cite{flandoli} and references therein for FPK equations on infinite dimensional state space). We may need finite dimensional approximations. We can project this parabolic PDE  according to either the $L^2$ direct metric (inducing the metric $\gamma(\theta)$) or, by deriving the analogous equation for $\sqrt{p_t}$, 
 $\frac{\partial \sqrt{p}}{\partial t} = \frac{1}{\sqrt{p}} {\cal L}_t^\ast\, p$, 
according to the Hellinger metric (inducing $g(\theta)$). 
The respective tangent linear projections
\begin{equation}\label{eq:projgamg} \frac{d}{dt} p(\cdot,\theta(t)) = \Pi^\gamma_{\theta(t)}\left[  {\cal L}_t^\ast\, p(\cdot,\theta)\right], \ \ 
\frac{d}{dt} \sqrt{p(\cdot,\theta(t))} = \Pi^g_{\theta(t)} \left[\frac{1}{\sqrt{p(\cdot,\theta)}} {\cal L}_t^\ast\, p(\cdot,\theta)\right] 
\end{equation}
 transform the PDE into a finite dimensional ODE for $\theta$ via the chain rule: 
\[ \frac{d}{dt} {\rho(p(\cdot,\theta_t))} =  
\sum_{j=1}^n
   \frac{\partial {\rho(p(\cdot,\theta))}}{\partial \theta_j} \dot{\theta}_j(t), \ \ \rho(p) = p \ \ \mbox{or} \ \ \rho(p) = \sqrt{p}. \]
The basic idea is illustrated in Figure \ref{fig:projections}.   
    
For brevity we write, for suitable functions $\varphi$, $E_\theta[\varphi] := \int \varphi(x) p(x,\theta) dx$.
The tangent projections in direct and Hellinger metric respectively yield, after an integration by parts using the fact that ${\cal L}^\ast$ is the formal adjoint  of $\cal L$:
%\begin{eqnarray*} 
%   \dot{ \theta}_t^i 
%   = \displaystyle 
%    \sum_{j=1}^m \gamma^{ij}(\theta_t)\; 
%   \int {{\cal L}_t^\ast\, p(x,\theta_t)}\; 
%   \frac{\partial p(x,\theta_t)}{\partial \theta_j} dx  , \ \ \theta^i_0 \ . \ \ 
%\end{eqnarray*}
%\begin{eqnarray*} 
%   \dot{\theta}_t^i 
%   = 
% \sum_{j=1}^m g^{ij}(\theta_t)\; 
%   \int \frac{{\cal L}_t^\ast\, p(x,\theta_t)}{p(x,\theta_t)}\; 
%   \frac{\partial p(x,\theta_t)}{\partial \theta_j}\; 
%   dx     , \ \ \theta^i_0 \ ,
%\end{eqnarray*}
%or, using integration by parts and the formal adjoint operator,
\begin{eqnarray}\label{eq:projboth} 
   \dot{ \theta}_t^i 
   = \displaystyle 
   \sum_{j=1}^n \gamma^{ij}(\theta_t)\; E_{\theta_t}
   \left[ 
   {\cal L}\left(\frac{\partial p(\cdot,\theta_t)}{\partial \theta_j}\right) \right], \     \dot{\theta}_t^i 
   =  \sum_{j=1}^n g^{ij}(\theta_t)\;     
E_{\theta_t}   \left[ 
   {\cal L}\left(\frac{\partial \log p(\cdot,\theta_t)}{\partial \theta_j}\right) \right]      
\end{eqnarray}
both equations starting at $\theta^i_0$.
\begin{remark}[Square roots and deformed logarithms]
Given the important role played by square roots in defining the Hellinger distance, one might have expected square roots to show up in the projected equation, namely in the second equation in \eqref{eq:projboth}. We emphasized the logarithm because this is particularly natural in view of our application to exponential families, but the square roots are still there. Indeed, an alternative and equivalent representation of the second equation in \eqref{eq:projboth} would be
\[  \dot{\theta}_t^i 
   =  \sum_{j=1}^n 4 \ g^{ij}(\theta_t)\; \int \left( \frac{{\cal L}_t^* p(\cdot,\theta_t)}{2 \sqrt{p(x,\theta_t)}} \frac{\partial \sqrt{p(x,\theta)}}{\partial \theta_j}\right) dx \]
Expressing the partial derivatives of square roots via the chain rule and integrating by parts gives immediately the second equation in \eqref{eq:projboth}. Finally, while here we use the two maps $p \mapsto \sqrt{p}$ and $p \mapsto p$, one might use different maps in the spirit of the theory of deformed logarithms, see \cite{naudts}, see also \cite{Newton1,Newton2,Newton3}.
\end{remark}

\subsection{Exponential families and expectation-to-canonical formula}
We now choose specific families to carry out the tangent projection. In particular, we project on the exponential family using the Hellinger distance and on the mixture family using the direct distance. 
The exponential families we consider are EF$(c)$, whose generic density is defined as $p(x,\theta):= \exp[\theta' c(x) - \psi(\theta)]$. The functions $c$ are the sufficient statistics of the family, the parameters $\theta \in \Theta \subset \Rx^n$ are the canonical parameters of the family. The quantity $\psi(\theta)$ is a normalizing constant needed for the density to integrate to one. 
Exponential families work well with the Hellinger/Fisher Rao choice because the tangent space has a simple structure: square roots do not complicate issues thanks to the exponential structure. 
A number of further results can be found or immediately derived  from Amari~\cite{amari85a} (Chapter 4) or 
Barndorff-Nielsen~\cite{barndorff-nielsen78a} (Theorem~8.1). Indeed, 
 the Fisher matrix has a simple structure: $\partial_{\theta_i,\theta_j}^2 \psi(\theta) = g_{ij}(\theta)$. The structure of the projection $\Pi^g$ is simple for exponential families. Finally, alternative coordinates, expectation parameters, defined via $\eta(\theta) = E_\theta [ c] = \nabla \psi =(\partial_{\theta_i} \psi(\theta))_i$ are available, with the two coordinate systems $\eta$ and $\theta$ being bi-orthogonal.  
Notice further that 
\begin{equation}\label{eq:detagdtheta} \partial_{\theta_i} \eta_j(\theta) = g_{i,j}(\theta), \ \   d \eta(\theta) = g(\theta) d\theta . 
\end{equation}
We further have $E_{\theta}\{c_i c_j\} 
   = \partial_{ij}^{2} \psi(\theta) + \partial_i \psi(\theta)   \partial_j \psi(\theta)$ and more generally
\begin{displaymath} 
   E_{\theta}\{c_{i_1} \cdots c_{i_k}\} 
   = \exp[ -\psi(\theta) ]\; 
   \frac{\partial^k  \exp[ \psi(\theta) ]}
   {\partial \theta_{i_1} \cdots \partial \theta_{i_k}}\; \ . 
\end{displaymath} 

To have a well-behaving matrix $g$ and good properties for the map $\theta \mapsto \psi(\theta)$ one typically requires that the sufficient statistics $(c_i)_i$ in the exponential family are linearly independent.

Consider now a special case for the family EF$(c)$. 
More specifically, we take the exponential polynomial manifold
EP$(n) := \{p(\cdot,\theta): \theta \in \Theta \subset \Rx^n\}$,
with $m$ an even positive integer and
with a linear combination of the monomials $x, x^2,\ldots,x^n$ 
in the exponent:
\begin{equation} \label{exfamily1}
p(x,\theta) = \exp\{\theta_1 x + ... + \theta_n x^n - \psi(\theta)\},
\ \ \theta_n < 0.
\end{equation}
In the PhD dissertation of Brigo (1996) \cite{BrigoPhD} (lemma 3.3.3), later partly published in Brigo and Hanzon (1998) \cite{brigoime}, the following result is introduced.  
\begin{theorem}[Expectation-canonical parameters algebraic relation for  EP$(n)$]\label{th:exptocan}
For the family EP$(n)$ with $n$ even positive integer, characterized by
$c_i(x) = x^i$, $i=1,\cdots,n$, $\theta_n < 0$, 
  the following recursion formula holds, with $\eta_0(\theta) := 1$.
   For any nonnegative integer $i$
\begin{eqnarray} \label{recexpar}
  && 
\eta_{n+i}(\theta) := E_{\theta}\{x^{n+i}\} 
  \\ \nonumber   
  && = - \frac{1}{n \theta_n} \left[ 
\begin{array}{ccccc} 
   (i+1) & \theta_1 & 2 \theta_2 & \cdots & (n-1)\theta_{n-1}
\end{array} \right]\; \left[ 
\begin{array}{c}
   \eta_i(\theta) \\ \eta_{i+1}(\theta) \\ \eta_{i+2}(\theta) \\ 
   \vdots \\ \eta_{i+n-1}(\theta) 
\end{array} \right]\ .  
\end{eqnarray}
   Moreover, the entries of the Fisher information matrix satisfy 
\begin{equation} \label{magic} 
   g_{ij}(\theta) 
   = \eta_{i+j}(\theta) - \eta_i(\theta)\, \eta_j(\theta)\ . 
\end{equation} 
Consequently, by defining the matrix
$M(\eta)$ as follows:
\begin{eqnarray}
   M_{i,j}(\eta) := \eta_{i+j}, \ \ \ i,j=1,2,\ldots,n
\end{eqnarray}
it is easy to verify that \eqref{recexpar} and the related lemma imply the following formula:
\begin{eqnarray} \label{theta_eta}
   \left[
   \begin{array}{c}
       \theta_1\\
       2 \theta_2 \\
       \vdots \\
       n \theta_n
   \end{array}
   \right]
   \ \ = \ \ -M(\eta)^{-1} \ \
   \left[
   \begin{array}{c}
       2 \eta_1\\
       3 \eta_2 \\
       \vdots \\
       (n+1) \eta_n
   \end{array}
   \right].
\end{eqnarray}
From this last equation it follows that we can recover 
{\em algebraically} the canonical parameters $\theta$ from the 
knowledge of the moments $\eta_1,\ldots,\eta_{2n}$ up to order $2n$.

\end{theorem}

{\sc Proof~:}
The recursion formula~(\ref{recexpar}) is obtained 
via integration by parts:
\begin{displaymath} 
\begin{array}{rcl} 

&&  \int_{-\infty}^{+\infty} x^i\; p(x,\theta)\,dx \\ \\
 && =  
   [\frac{x^{i+1}}{i+1}\; p(x,\theta)]_{-\infty}^{+\infty} 
   - \int_{-\infty}^{+\infty} \frac{x^{i+1}}{i+1}\; 
   [\theta_1 + 2 \theta_2\, x + \cdots + n \theta_n\, x^{n-1}]\; 
   p(x,\theta)\,dx \\ \\ 
 &&  = \displaystyle 
   0 - \frac{1}{i+1}\; 
   E_{p(\cdot,\theta)}\{\theta_1\, x^{i+1} + 2 \theta_2\, x^{i+2} 
   + \cdots + n \theta_n\, x^{i+n} \}\ , 
\end{array} 
\end{displaymath} 
from which the formula and the other results follow easily. \qed

The above results for EP$(n)$ solve the problem of recovering the density and the canonical parameters $\theta$ from knowledge of the expectation parameters $\eta$.  The opposite direction is straightforward: from (\ref{exfamily1}) it is clear that the canonical
parameters $\theta$ permit to express the densities of $EP(n)$ explicitly. 

For a study of such  procedure, in a slightly different context, and for a 
comparison with several alternatives, including a Newton method,
see Borwein and Huang (1995)  \cite{BorweinHuang95}.
Further investigations into this so called polynomial moment
problem are called for. Better insight into the geometry of the
manifolds $EP(n)$ is likely to be helpful, especially to understand
the behaviour of the various algorithms at the boundary of the
manifold where $\theta_n$ is close to zero. 

This concludes our summary of exponential families results.

\subsection{Simple mixture families}
In dealing with exponential families, we used the Hellinger metric $d_H$ since this resulted in a number of good properties summarized in the previous section. For the direct metric $d_D$, instead, it may be more interesting to  project on simple mixture families. We define a \emph{simple mixture family} as follows. Given $n+1$ fixed squared integrable probability densities  $\uq = [q_1,q_2,\ldots,q_{n+1}]'$, define $\hat{\theta}(\theta) := [ \theta_1, \theta_2, \ldots, \theta_n, 1 - \theta_1- \theta_2 - \ldots - \theta_n]'$ for all $\theta \in \mathbb{R}^n$, and set $q = [q_1,q_2,\ldots,q_{n}]'$.  The simple mixture family (on the simplex) is defined as 
\begin{equation} \mbox{SM}(\uq) := \{p(\cdot,\theta) = \hat{\theta}(\theta)' \uq(\cdot), \ \ \theta \in \Theta \subset \Rx^n\}, \ \ 
\end{equation}
\[\Theta\ \  \mbox{open subset of}\  \  \{\theta \in \Rx^n: \theta_i \in [0,1]\ \mbox{for all} \ \ i \ \ \mbox{and} \ \ 0<  \sum_{i=1}^n \theta_i < 1\}.\]

If we consider the $L^2$ based $d_D$ distance with metric $\gamma(\theta)$, the metric $\gamma(\theta)$ itself and the related projection become very simple and do not depend on the point $\theta$. For example, $\frac{\partial p(\cdot,\theta)}{\partial \theta_i} = q_i - q_{n+1}$, $\gamma_{i,j} = \langle q_i-q_{n+1},q_j - q_{n+1} \rangle$. 
% \ \ \mbox{and} \ \  \gamma_{ij}(\theta) = \int (q_i(x) - q_{m+1}(x))(q_j(x) - q_{m+1}(x)) dx \ \ \ \  \]
%requiring no inline numeric integration. 
Accordingly, the tangent space at $p(\cdot,\theta)$ does not depend on $\theta$ and is given by
$\mbox{span}\{ q_1-q_{n+1}, q_2-q_{n+1}, \cdots, q_n-q_{n+1}\}$. 

We now introduce expectation parameters for the simple mixture family SM$(\uq)$.
Define
\[ m(\theta) = E_\theta[q- q_{n+1} \uni ]\] 
where $\uni$ is a vector with $m$ components equal to $1$. The reason for subtracting $q_{n+1}$ will be clear when we consider the metric projection later. 
A quick calculation shows that
\begin{equation}\label{dmgammadtheta} m(\theta) = \gamma \theta + \beta, \ \ \beta_i = E_{q_{n+1}}[ q_i - q_{n+1}], \ \ i=1,\ldots,n. 
\end{equation}
Note that neither $\gamma$ nor $\beta$ depend on $\theta$, so that
\begin{equation}\label{dmgammadtheta2} 
d m(\theta) =  \gamma\  d\theta. 
\end{equation}
This is the simple mixture counterpart of the exponential families relationship \eqref{eq:detagdtheta}. 
 
Since $m$ is an alternative parameterization for the mixture family, we will denote the simple mixture density characterized by $m$ with $p(\cdot;m)$. 

%Also the $L^2$ projection becomes particularly simple and is equivalent both to a classic Galerkin method and to an assumed density approximation, as we will show below.

\subsection{Vector-field tangent projection on mixture \& exponential families}\label{sec:vectorfieldprojhellinger}

Starting from the simple mixture family with the direct metric, we now specialize our tangent-space projection equations to $(SM(\uq),d_D)$.  The first Eq. (\ref{eq:projboth}) specializes to 
\begin{eqnarray}\label{eq:projl2mix}
   \dot{ \theta}_t^i 
   = \displaystyle 
   \sum_{j=1}^n \gamma^{ij}\; \sum_{k=1}^{n+1}
   \hat{\theta}_k(t) \int  
   ({\cal L}\left(q_j(x)-q_{m+1}(x)\right)) \  (q_k(x)-q_{m+1}(x) ) dx  , \ \  \theta^i_0 \ , \ \ 
\end{eqnarray}
 which is a linear equation, see \cite{armstrongbrigomcss} for the more general case of the filtering problem SPDE. The quantity $\theta^i_0$ in the right hand side denotes the given initial condition. See also \cite{brigopistone} for more details, and in particular see \cite{santacroce} for mixture and exponential models in Orlicz spaces. See also Ahmed \cite{ahmedbook}, Chapter 14, Sections 14.3 and 14.4 for a summary of the Galerkin method for the Zakai equation, of which FPK is a special case, and see also \cite{ahmed97}.
 
For the derivation of Eq. \eqref{eq:projl2mix} to make sense we need to ensure that we are indeed projecting a $L^2$ vector field from the related FPK equation. 
 This is ensured if we require that ${\cal L}_t^* p(\cdot,\theta)$ can be tangent-linearly projected via $\Pi^\gamma$ as an $L^2$ vector for $p(\cdot,\theta) \in $SM$(\uq)$. 
 %For this to happen, we require
\begin{equation}\label{eq:L2existdirect} \sup_{t \ge 0} \|  {\cal L}_t^\ast p(\cdot,\theta)   \|_2 < \infty \ \ \mbox{for all} \ \theta \in \Theta.
\end{equation}
Given that ${\cal L}_t^* p(\cdot,\theta) =  \hat{\theta}' {\cal L}_t^* \uq$ (where ${\cal L}_t^*$ is meant to be applied component-wise when applied to vectors), it is enough to ensure that each ${\cal L}_t^* q_i$ can be projected via the $L^2$ inner product. By inspection of 
${\cal L}_t^* q_i$ we can formulate the following
\begin{proposition}\label{rem:L2existmix}(Sufficient conditions for the direct $L^2$ structure to apply to the mixture case). Sufficient conditions guaranteeing the $L^2$  condition \eqref{eq:L2existdirect} 
in SM$(\uq)$ are the following: $f_t$ and its first derivatives, $a_t$ and its first and second derivatives, and $\uq$ with their first and second derivatives have at most polynomial growth, and densities $\uq$ together with their first and second derivatives integrate any polynomial.  
\end{proposition}

An immediate example of simple mixtures satisfying the above polynomial and integrability conditions are Gaussian mixtures.

%
%
%
%

%\medskip

%\subsection{$g(\theta)$-vector field projection on exponential families \& MLE}
Moving to exponential families, we will project in Hellinger distance/Fisher-Rao metric, as hinted above. The tangent linear projection of the FPK equation in Fisher metric has been introduced first in \cite{brigogyor}. 
%First of all we refer to the the Hellinger projection above, the second Eq. in \eqref{eq:projboth}, 
%as to the Hellinger {\emph{``vector field projection''}}. 
%This is because we project the $L^2$ vector field of the Fokker Planck equation for $\sqrt{p}$ instant by instant onto the tangent space of EF$(c)$, thus obtaining a vector field in EF$(c)$. 
For the second Eq. in \eqref{eq:projgamg} to hold we need to ensure that $({\cal L}_t^\ast p(\cdot,\theta) )/\sqrt{p(\cdot,\theta)}$ is indeed an $L^2$ vector for all $t$ and $\theta$. This holds in turn if the condition
\begin{equation}\label{eq:L2exist} \sup_{t \ge 0} E_\theta \left[ \left| \frac{{\cal L}_t^\ast p(\cdot,\theta) }{p(\cdot,\theta)} \right|^2\right] < \infty \ \ \mbox{for all} \ \theta \in \Theta
\end{equation}
holds. Again by inspection, we have the following

\begin{proposition}\label{rem:L2exist}(Sufficient conditions for the Hellinger $L^2$ structure to apply). Sufficient conditions guaranteeing condition \eqref{eq:L2exist} in EF$(c)$ are the following: $f_t$ and its first derivatives, $a_t$ and its first and second derivatives, and $c$ with its first and second derivatives have at most polynomial growth, and densities in EF$(c)$ integrate any polynomial.  
\end{proposition}

See also \cite{brigo99} for a proof and a detailed discussion, including   conditions under which all vector fields are well defined and the tangent projection is well defined, see in particular Theorem 5.4 in \cite{brigo99} in the special case $h=0$.

\begin{eqnarray}\label{as:exist} \mbox{{\bf Assumption}. From now on, we assume sufficient conditions above for}\ \eqref{eq:L2existdirect}\\ \nonumber \mbox{in the mixture case and for}\ \eqref{eq:L2exist}\ \mbox{in the exponential case hold.
For example,} \\ \nonumber \mbox{  we might assume conditions  given in Propositions} \  \ref{rem:L2existmix} \ \mbox{or} \ \ref{rem:L2exist}\ \mbox{to hold.}
\end{eqnarray}

Again for brevity, we write $E_\eta[\varphi] := \int \varphi(x) p(x;\eta(\theta)) dx$. 
When projecting on EF$(c)$ the second Eq. (\ref{eq:projboth}) specializes into
\begin{equation}\label{FPEPRO:single}
   \dot{\theta}_t^i 
   =  \sum_{j=1}^m g^{ij}(\theta_t)\; 
   E_{\theta_t} 
   \left[ 
   {\cal L} c_i  \right]  , \ \theta^i_0 \ , \ \ \  \mbox{or} \ \ \ 
   \dot{\eta}_t^i 
   =   
   E_{\eta_t} 
   \left[ 
   {\cal L} c_i  \right]   ,  \ \eta^i_0 \ 
\end{equation}
where the second equation has been obtained from the first by recalling that  $d \eta(\theta) = g(\theta) d\theta$.

We know that the tangent-linear orthogonal projection is giving us the locally optimal approximation for the chosen metric. A natural question when projecting is how good the projection is locally or, in other terms, how well does the chosen finite dimensional family approximate the infinite dimensional evolution locally? Indeed, one would like to have a measure for how far the projected evolution is, locally, from the original one. We now define a local projection residual as the $L^2$ norm of the FPK infinite dimensional vector field minus its finite-dimensional orthogonal projection. Define the vector field minus its tangent projection and the related norm as
\[ \varepsilon_t(\theta) := \frac{{\cal L}_t^\ast p(\cdot,\theta)}
          { 2 \sqrt{p(\cdot,\theta)}} - \Pi^g_\theta
  \left[\frac{{\cal L}_t^\ast p(\cdot,\theta)}
          {2 \sqrt{p(\cdot,\theta)}}\right], \ \ R^2_t : = {\| \varepsilon_t(\theta) \|_2}^2  \]
The projection residual $R_t$ can be computed jointly with the projected equation evolution \eqref{FPEPRO:single} to have a local measure of the goodness of the approximation involved in the projection.

Monitoring the projection residual and its peaks can be helpful in tracking the local projection method performance, see also  \cite{brigo99} for examples of $L^2$-based projection residuals in the more complex case of the Kushner-Stratonovich equations of nonlinear filtering. However, the projection residual only allows for a local approximation error numerical analysis. 
To illustrate this, assume for a moment that time is discrete $1,2,3,\ldots$ and consider again Figure \ref{fig:projections}. To make the point, we are artificially separating tangent projection and propagation and the local and global errors. This is not completely precise but allows us to make an important point on our method. If we start from the manifold with a $p_0 = p(\cdot,\theta^0) \in$ EF$(c)$ (omitting square roots in the notation and with upper indices denoting time), the FPK vector field driven by ${\cal L}^\ast$ will move us out of the manifold as the $L^2$ vector related to ${\cal L}^\ast p^0$  will not be tangent to EF$(c)$. We then project this $L^2$ vector on the tangent space of EF$(c)$ and follow it, obtaining a new $p(\cdot,\theta^1)$ in the manifold. Now we continue, again the FPK vector field driven by ${\cal L}^\ast p(\cdot,\theta^1)$ would bring us out of EF$(c)$, and to avoid this we project it onto the tangent space and follow the projected vector, obtaining $p(\cdot,\theta^2)$. The crucial point here is that this second step was done starting from an approximate point $p(\cdot,\theta^1)$ rather than from the true FPK density $p_1$. This means that, besides the local projection error measured by $R_t$, we have a second error coming from the fact that we start the projection from the wrong point. 
If we leave the global approximation error analysis aside for a minute, the big advantage of the above method is that it does not require us to know the true solution of the FPK equation to be implemented. Indeed, Equation \eqref{FPEPRO:single} works perfectly well without knowing the true solution $p_t$.

\section{Globally optimal approximation to FPK}
We now investigate whether it is possible to say something on the globally optimal approximation of the solution of the FPK equation on a family $\{p(\cdot, \theta), \ \theta \in \Theta\}$. This will turn out to be related to a form of moments matching, or matching of expectations. In trying to find related results, we will find a local approximation based on a localized version of the globally optimal approximation via the assumed density approximation, and we will find also a tractable case for the global approximation based on the FPK diffusion operator eigenfunctions.

\subsection{Metric projection}

\begin{figure}[htp]
\begin{centering}
\includegraphics[scale=0.5]{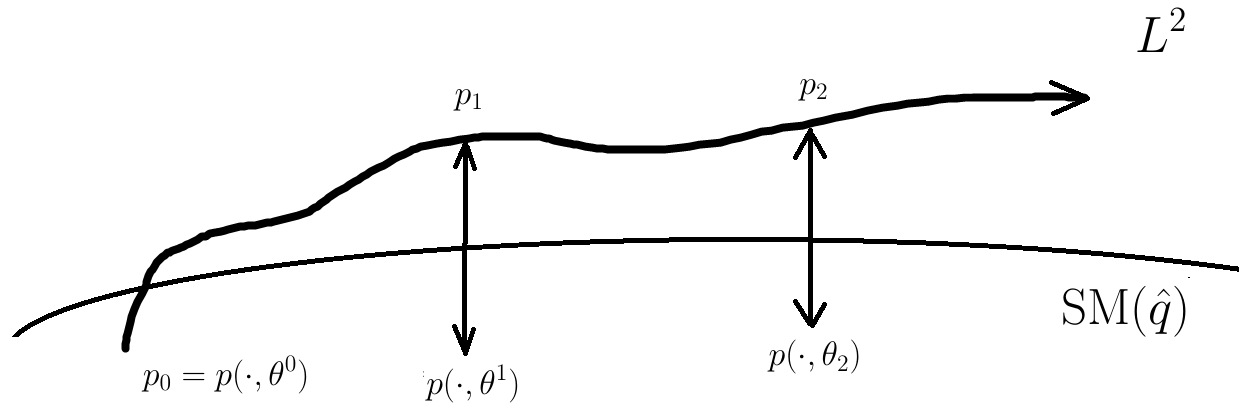}
\caption{Metric projection of the FPK Equation on the manifold SM$(\uq)$ in direct metric. We start from an initial condition $p_0$ in the manifold. $p_0,p_1,p_2$ describes the exact solution in $L^2$, whereas $p_0,p(\cdot,\theta^1), p(\cdot,\theta^2)$ describes the evolution resulting from the metric projection. Notice the difference with Figure \ref{fig:projections}, where we project vectors rather than minimizing the distance from the manifold.} 
\label{fig:mprojections}
\end{centering}
\end{figure}

Now, to study the global error, we introduce a second projection  method. We call this method
{\emph{``metric projection''}}, since here we will project directly the $L^2$ densities or their square roots  onto SM$(\uq)$ or EF$(c)$, rather than projecting these densities evolutions. In particular, we will make no use of tangent spaces. In geometry this is called indeed ``metric projection'', as opposed to the tangent space linear projection we used so far. This is illustrated in Figure \ref{fig:mprojections}.

Metric projection will require us to know the true solution, so as an approximation method it will be pointless. However, it will help us with the global error analysis, and a modification of the method based on the assumed density approximation will allow us to find an algorithm for a local approximation that does not require the true solution. 

\subsection{Metric projection on exponential families in relative entropy}

The metric projection method for exponential families works as follows. Starting again from the manifold with a $p_0 = p(\cdot,\theta^0) \in$ EF$(c)$, the FPK vector field driven by ${\cal L}^\ast$ will move us out of the manifold; we follow this vector field and reach $p_1$.  
To go back to EF$(c)$ we project $p_1$ onto the exponential family by minimizing the relative entropy, or Kullback Leibler information of $p_1$ with respect to EF$(c)$, finding the orthogonal projection of $p_1$ on the chosen manifold. It is well known that the orthogonal projection in relative entropy is obtained by matching the sufficient statistics expectations of the true density. Namely, the projection is the particular exponential density $p(\cdot;\eta^1) \in EF(c)$ with $c$-expectations $\eta^1 = \teta^1$ where $\teta^1 = E_{p_1}[c]$ are the $c$-expectations of the density $p_1$ to be approximated. 
          See for example  Kagan et al. (1973) \cite{kagan}, Theorem 13.2.1, or \cite{brigoime} for a quick proof and an application to filtering in discrete time. 
For convenience, we briefly prove the theorem below.

\begin{theorem}[Moment matching \& entropy minimization for exponential families]\label{th:mmexp}
Suppose we are given the exponential family EF$(c)$ and a $L^1$ probability density $p$ outside EF$(c)$. Suppose the matrix $g(\theta)$ for EF$(c)$ is positive definite. Then the minimization problem
\[  \min_{\theta \in \Theta} K(p,p(\cdot,\theta)) \]
that consists in finding the density in $p(\cdot;\eta) \in EF(c)$ that is closest, in relative entropy, to $p$ has a unique solution characterized by the ``moment matching'' or ``$c$-expectations matching'' conditions
\[ \eta_i = \teta_i \ \ \mbox{where} \ \ \teta_i= E_p[c_i], \  \eta_i = E_\theta[c_i]  \ \ \mbox{for} \ \ i=1,\ldots,n .\] 
In other words, the EF$(c)$ density that is closest to $p$ in relative entropy is the density in EF$(c)$ that shares the $c$-expectations (or $c$-moments) with $p$.  
\end{theorem}          
The proof is immediate. Write
\[ K(p,p(\cdot,\theta)) = \int (\ln p(x) - \theta' c(x) + \psi(\theta)) p(x) dx \]
and note that the gradient with respect to $\theta$ reads
\[ \frac{\partial K(p,p(\cdot,\theta))}{\partial \theta_i} = \int ( - c_i(x) + \partial_{\theta_i} \psi(\theta)) p(x) dx = - E_p[c_i] + \eta_i(\theta)  \]
where we used the basic property of exponential families introduced above, namely $\partial_{\theta_i} \psi = \eta_i = E_\theta[c_i]$. A necessary condition for minimality is given by setting the gradient to zero. This yields
\[ E_p[c_i] = \eta_i(\theta^*) \ \ \mbox{for} \ \ i=1,\ldots, n.\]
To check that the condition is also sufficient we need to show that the Hessian is positive definite. Compute
\[ \frac{\partial^2 K(p,p(\cdot,\theta))}{\partial \theta_i \partial \theta_j} = \int (  \partial^2_{\theta_i,\theta_j} \psi(\theta)) p(x) dx =    \partial^2_{\theta_i,\theta_j} \psi(\theta) = g_{i,j}(\theta) \]
again by the basic properties of exponential families we saw above. 
%Since $g$ is the Hessian of the strictly convex function $\theta \mapsto \psi(\theta)$ it is positive definite. 
Since we are assuming $g$ to be positive definite, this concludes the proof. 
\qed

          As we hinted above, we know that EF$(c)$, besides $\theta$, admits another important coordinate system, the expectation parameters $\eta$. If one defines 
          $\eta(\theta) = E_{p(\theta)} [ c]$ as above, 
          then $d \eta(\theta) = g(\theta) d\theta$ where $g$ is the Fisher metric, as we have seen earlier.
          Thus, we can take the $\eta^1$ above coming from the true density $p_1$ and look for the exponential density $p(\cdot;\eta^1)$ sharing these $c$-expectations. This will be the closest in relative entropy to the true $p_1$ in EF$(c)$.
We then continue moving forward in time, iterating this algorithm. 

The advantage of this method compared to the previous vector field based one is that we find at every time the best possible approximation (``maximum likelihood'') of the true solution in EF$(c)$. The disadvantage is that in order to compute the projection at every time, such as for example $\eta^1 = E_{p_1}[c]$, we need to know the exact solution $p_1$ at that time.  
However,  it turns out that we can get back an algorithm that does not depend on the exact solution if we invoke the assumed density approximation. 

\subsection{Back to local: assumed density approximation for exponential families}

 This works as follows. 
 Differentiate both sides of $\teta_t = E_{p_t}[c]$ to obtain
 \[ \frac{d}{dt}{\teta}_t = \frac{d}{dt} \int c(x) p_t(x) dx = \int c(x) \frac{\partial p_t(x)}{\partial t} dx = \int c(x) {\cal L}^\ast_t p_t(x) dx  = E_{p_t} [ {\cal L} c] \] 
so that 
\begin{equation}\label{eq:etatrue}  \frac{d}{dt}{\teta}_t = E_{p_t} [ {\cal L} c]. 
\end{equation} 
This last equation is not a closed equation, since $p_t$ in the right hand side is not characterized by $\teta_t$. Thus, to be solved this equation should be coupled with the original FPK equation for $p_t$ and we would still be in infinite dimension. However, at this point we can close the equation by invoking the assumed density approximation (ADA). Implement the following replacements in Equation \eqref{eq:etatrue}. 
\[ \mbox{ADA: \ in Eq. \eqref{eq:etatrue}} \ \  \teta_t \rightarrow \eta_t \ \mbox{(left hand side)}, \ \  \   p_t \rightarrow p(\cdot;\eta_t) \ \mbox{(right hand side)}. \]
We obtain
\begin{equation}\label{eq:etaode} {\dot{\eta}}_t = E_{\eta_t} [ {\cal L} c].
\end{equation}
 This is now a finite dimensional ODE for the expectation parameters. It does not require the true solution to be implemented but the arbitrary replacement implies that we have compromised global optimality. However, perhaps surprisingly, we still have local optimality in the same sense as we had with the tangent space projection. Indeed, this last equation is the same as our earlier vector field based projected equation
 \eqref{FPEPRO:single}. This result had been proven for nonlinear filtering in \cite{brigo99}. Intuitively, the result is related to the fact that the Fisher-Rao metric and relative entropy are infinitesimally equivalent, see Eq \eqref{eq:klig}.
 \begin{theorem}[In EF$(c)$ Hellinger-Fisher-Rao tangent projection = ADA] Assumption 
\eqref{as:exist} in force. 
Closing the evolution equation for the relative entropy point projection of the FPK solution onto EF$(c)$ by forcing an exponential density on the right hand side is equivalent to the locally optimal Hellinger approximation based on the vector field tangent linear projection in Fisher metric. 
 \end{theorem} 

\subsection{Global optimality for exponential families: ML \& $\cal L$ eigenfunctions}

We can now attempt an analysis of the error between the best possible projection $\teta_t$ and the vector field based (or equivalently assumed density approximation based) projection ${\eta}$. To do this, write
$\epsilon_t := \teta_t - \eta_t$, 
expressing the difference between the best possible approximation and the vector field tangent projection / assumed density one, in expectation coordinates. 
Differentiating we see easily that
$\dot{\epsilon}_t =  (E_{p_t} [ {\cal L} c] -  E_{p({\eta}_t)} [ {\cal L} c] )$. 
Now suppose that the $c$ statistics in EF$(c)$ are chosen among the eigenfunctions of the operator ${\cal L}$, so that
${\cal L} c =  - \Lambda c$, 
where $\Lambda$ is a $n \times n$ diagonal matrix with the eigenvalues corresponding to the chosen eigenfunctions. 
Substituting, we obtain
\[ \dot{\epsilon}_t =  - \Lambda (E_{p_t} [  c] -  E_{p({\eta}_t)} [  c] )  \ \ \ \mbox{or} \ \ \  \dot{\epsilon}_t =  - \Lambda  \epsilon_t \Rightarrow  \epsilon_t = \exp(- \Lambda t) \epsilon_0  \]
so that if we start from the manifold ($\epsilon_0=0$) the error is always zero, meaning that the vector field tangent projection gives us the best possible Maximum Likelihood (ML) approximation. If we don't start from the manifold, ie if $p_0$ is outside EF$(c)$, then the difference between the vector field approach and the best possible approximation dies out exponentially fast in time provided we have negative eigenvalues for the chosen eigenfunctions. This leads to the following 
\begin{theorem}[Global optimality/ML for FPK  \& Fisher-Rao projection.]
Consider the Fokker-Planck-Kolmogorov equation and an exponential family EF$(c)$.
Assumption \eqref{as:exist} in force. 
The vector field tangent linear projection approach leading to the locally optimal approximation \eqref{FPEPRO:single} in Hellinger distance in EF$(c)$ provides also the global optimal approximation of the Fokker-Planck-Kolmogorov equation solution in relative entropy in the family EF$(c)$, provided that the sufficient statistics $c$ are chosen among the eigenfunctions of the adjoint operator $\cal L$ of the original Fokker-Planck-Kolmogorov equation, and provided that EF$(c)$ is an exponential family when using such eigenfunctions. In other words, under such conditions the Fisher Rao vector field projected equation \eqref{FPEPRO:single} provides the exact maximum likelihood density for the solution of the Fokker-Planck-Kolmogorov equation in the related exponential family. 
\end{theorem}

Before briefly discussing the eigenfunctions result, we consider an analogous setting and derivation for mixture families

\subsection{Metric projection on mixture families in $L^2$ direct distance}
We begin by deriving the metric projection. Our result is the analogous of Theorem \ref{th:mmexp} for mixture families and direct metric.

\begin{theorem}[Moment matching \& $L^2$ distance minimization for simple mixtures]\label{th:mmmix}
Suppose we are given the simple mixture family SM$(\uq)$  and a $L^2$ probability density $p$ outside SM$(\uq)$. Suppose the matrix $\gamma$ for the family SM$(\uq)$ is positive definite. Then the minimization problem
\[  \min_{\theta \in \Theta} \| p - p(\cdot,\theta) \|_2 \]
that consists in finding the density in SM$(\uq)$ that is closest in $L^2$ direct distance $d_D$ to $p$ has a unique solution characterized by the ``$q_i-q_{n+1}$ moment or expectation matching'' conditions
\[\tim = m \] where 
\[ \tim = E_p \left[ q - \uni q_{n+1} \right] \ \mbox{and} \ \  m = E_\theta \left[ q - \uni q_{n+1} \right] = \gamma \theta + \beta,  
 \ \ \beta = E_{q_{n+1}}[ q - q_{n+1} \uni] ,\]
provided that the solution $\theta^* = \gamma^{-1}(m-\beta) \in \Theta$.  
\end{theorem}          
The proof of the theorem is immediate and similar to the the proof given in Theorem \ref{th:mmexp}. One minimizes
\[ \int (p(x) -\hat{\theta}'(\theta) \uq(x) )^2 dx \]
by taking partial derivatives with respect to $\theta$'s and setting them to zero. The Hessian matrix is $\gamma$ and is by assumption positive definite. 

Similarly to what we have seen for exponential families, we point out that the metric projection solution $m$ is difficult to obtain. Compared to the vector field tangent projection method, the advantage is again that we find at every time the best possible approximation of the true solution in SM$(\uq)$. The disadvantage is that in order to compute the metric projection at every time, such as for example $m^1 = E_{p_1}[q-q_{m+1}\uni]$, we need to know the exact solution $p_1$ at that time.  
However, also in this case it turns out that we can get back an algorithm that does not depend on the exact solution if we invoke the assumed density approximation. 

\subsection{Back to local: assumed density approximation for mixture families}

 This works as follows. The metric projection of the FPK equation in $d_D$ metric on SM$(\uq)$ is defined by
 \[  \tim_t = E_{p_t}[q-q_{n+1}\uni] .\]
 Differentiate both sides to obtain
\[ \frac{d}{dt} \tim_t = \frac{d}{dt} \int [q(x)-q_{n+1}(x)\uni] p_t(x) dx
= \int [q(x)-q_{n+1}(x)\uni] ({\cal L}^*_t p_t)(x) dx = \]\[=
\int ({\cal L}(q-q_{n+1}\uni))(x) p_t(x) dx\]
where we substituted the FPK equation and we used integration by parts. 
  We have obtained
\begin{equation}\label{eq:timtrue}
 \frac{d}{dt} \tim_t = E_{p_t}\left[{\cal L}(q-q_{n+1}\uni)\right] .
\end{equation}
This equation is not very helpful as an approximation, since we need to know the exact $p_t$ to compute the right hand side. However, we can apply the ADA in this case too.
\[ \mbox{ADA: \ in Eq. \eqref{eq:timtrue}} \ \  \tim_t \rightarrow m_t \ \mbox{(left hand side)}, \ \  \   p_t \rightarrow p(\cdot;m_t) \ \mbox{(right hand side)}. \]
We obtain
\begin{equation}\label{eq:etaode} {\dot{m}}_t = E_{m_t} [ {\cal L} (q-q_{n+1}\uni)].
\end{equation}
 This is a finite dimensional ODE for the expectation parameters. Also in the mixture case, perhaps surprisingly, we still have local optimality in the same sense as we had with the tangent space projection. Indeed, this last equation is the same as our earlier vector field based projected equation \eqref{eq:projl2mix}, as one notices immediately by applying 
$d m(\theta) = \gamma \ d \theta$.  
 Intuitively, the result is related to the fact that we are using the same $d_D$ metric both in the vector field projection and in the metric projection.
 \begin{theorem}[In SM$(\uq)$ we have $d_D$ tangent projection = ADA] Assumption 
\eqref{as:exist} in force. 
Closing the evolution equation for the $d_D$ metric projection of the Fokker-Planck-Kolmogorov solution onto SM$(\uq)$ by forcing a mixture density on the right hand side is equivalent to the locally optimal $d_D$ approximation based on the vector field tangent linear projection in $L^2$ metric.  
 \end{theorem}

\subsection{Global optimality for mixture families: $\cal L$ eigenfunctions and Galerkin methods}

We can now attempt an analysis of the error between the best possible projection $\tim_t$ and the vector field based (or equivalently assumed density approximation based) projection $m_t$. To do this, write
$\varepsilon_t := \tim_t - m_t$, 
expressing the difference between the best possible approximation and the vector field tangent projection / assumed density one, in expectation coordinates. 
Differentiating we see easily that
$\dot{\varepsilon}_t =  (E_{p_t} [ {\cal L}(q-q_{n+1}\uni) ] -  E_{m_t} [{\cal L} (q-q_{n+1}\uni))$. 
Now suppose that the $q_i-q_{n+1}$ mixture components of SM$(\uq)$ are chosen among the eigenfunctions of the operator ${\cal L}$, so that
${\cal L} (q -q_{n+1}\uni) =  - \Lambda (q -q_{n+1}\uni)$, 
where $\Lambda$ is a $n \times n$ diagonal matrix with the eigenvalues corresponding to the chosen eigenfunctions. 
Substituing, we obtain
\[ \dot{\epsilon}_t =  - \Lambda (E_{p_t} [ q -q_{n+1}\uni] -  E_{m_t} [q -q_{n+1}\uni ] )  \ \ \ \mbox{or} \ \ \  \dot{\epsilon}_t =  - \Lambda  \epsilon_t \Rightarrow  \epsilon_t = \exp(- \Lambda t) \epsilon_0  \]
so that, as in the exponential case, if we start from the manifold ($\epsilon_0=0$) the error is always zero, meaning that the vector field tangent projection gives us the best possible $d_D$ metric projection approximation. If we don't start from the manifold, ie if $p_0$ is outside SM$(\uq)$, then the difference between the vector field approach and the best possible approximation dies out exponentially fast in time provided we have negative eigenvalues for the chosen eigenfunctions. This leads to the following

\begin{theorem}[Global optimality for FPK  \& $d_D$ tangent projection]
Consider the Fokker-Planck-Kolmogorov equation and a simple mixture  family SM$(\uq)$.
Assumption \eqref{eq:L2existdirect} in force. 
The vector field tangent linear projection approach leading to the locally optimal approximation \eqref{eq:projl2mix} in $d_D$ direct $L^2$ metric in SM$(\uq)$ provides also the global optimal approximation of the FPK equation solution in $d_D$ in the same family, provided that the mixture components $q-q_{n+1}\uni$  are chosen among the eigenfunctions of the adjoint operator $\cal L$ of the original FPK equation, and provided that SM$(\uq)$ is a simple mixture family when using such eigenfunctions. 
%
%  h e r e 
%
%
\end{theorem}

Finally, in \cite{armstrongbrigomcss} we have shown that the vector field projection of the filtering SPDE in direct metric on simple mixtures is equivalent to a Galerkin method with basis functions given by the mixture components. Given that the FPK equation is a special case of the filtering equation, we can immediately derive the analogous result.

The Galerkin approximation is derived by approximating the exact solution $p_t$ of the FPK equation with a linear combination of basis functions $\phi_i(x)$, namely
\begin{equation}\label{eq:galerkin:approx} \tilde{p}_t(x) := \sum_{i=1}^{\ell} c_i(t) \phi_i(x) .
\end{equation}
The idea is that the $\phi_i$ can be extended to indices $\ell+1,\ell+2,\ldots,+\infty$ to form a basis of $L^2$. 

Now rewrite the FPK equation in weak form, using test functions:  
\[   \langle  -  \partial_t {p}_t + {\cal L}_t^\ast\, {p}_t\ ,  \xi \rangle = 0 \]
for all smooth $L^2$ test functions $\xi$ such that the inner product exists.

We replace this equation  with the equation
\[   \langle  -  \partial_t \tilde{p}_t + {\cal L}_t^\ast\, \tilde{p}_t
   , \ \phi_j \rangle = 0, \ \ j=1,\ldots,\ell . \]

By substituting Equation (\ref{eq:galerkin:approx}) in this last equation, using the linearity of the inner product in each argument and by using integration by parts we obtain easily an equation for the coefficients $c$, namely
\begin{eqnarray}\label{eq:galerkinmix} \sum_{i=1}^\ell \langle \phi_i,\phi_j \rangle \dot{ c}_i =   
\sum_{i=1}^\ell \langle \phi_i, {\cal L} \phi_j \rangle c_i   
\end{eqnarray}
We see immediately by inspection that this equation  coincides with the vector field projection \eqref{eq:projl2mix} of the FPK equation  if we take 
\[\ell = n+1, \  c_i = \theta_i   \ \mbox{and} \  \phi_i = q_i - q_{n+1}   \ \mbox{for} \ i=1,\ldots,n, \ \mbox{and} \  c_{n+1} = 1, \ \phi_{n+1} = q_{n+1}.\] 

We have thus proven the following 
\begin{theorem}
For simple mixture families SM$(\uq)$, the $d_D$ vector field projection approximation (\ref{eq:projl2mix}) coincides with a Galerkin method (\ref{eq:galerkinmix}) where the basis functions are the mixture components~$q$.
\end{theorem}

\section{Conclusions and further work}
We presented several methods to approximate the infinite dimensional solution of the Fokker-Planck-Kolmogorov equation with a finite dimensional probability distribution. The approximations we proposed are either locally or globally optimal, and consist of the vector field projection approximation and of the metric projection approximation. The two turn out to be related in a way that is similar to how the tangent space projection is related to the metric projection in general. Indeed, as the tangent projection is a linearization for small distances of the metric projection, so the vector field projection turns out to be equivalent to a localization of the metric projection based on the assumed density approximation. After exploring this setting for exponential and mixture families and further clarifying their relationship with maximum likelihood and Galerkin methods, we found a special case where one can implement the globally optimal approximation and this turns out to coincide with the vector field projection. This special case is based on the eigenfunctions of the original equation. In this case one obtains the exact maximum likelihood density.    
The choice of the specific exponential or mixture family, and the choice or availability of suitable eigenfunctions in particular is not always straightforward and requires further work. In this paper we focused on clarifying the relationship between different methods and on deriving the eigenfunctions result. For an initial discussion on eigenfunctions of the FPK equation see \cite{pavliotis} and \cite{liberzon} for the case of the Fokker-Planck-Kolmogorov equation for a linear SDE. For example, in the one dimensional case $N=1$ where the diffusion is on a bounded domain $[\ell, r]$ with reflecting boundaries and strictly positive diffusion coefficient $\sigma$ then the spectrum of the operator ${\cal L}$ is discrete, there is a stationary density and eigenfunctions can be expressed with respect to this stationary density, that could be taken as background measure, see \cite{brigopistone}.   
For the case $N>1$ special types of FPK equations allow for a specific eigenfunctions/eigenvalue analysis, see again \cite{pavliotis}. Further research is needed to explore the eigenfunctions approach in connection with maximum likelihood. 
In particular, the present paper might have some links with the variational approach to the Fokker-Planck-Kolmogorov equation started in \cite{otto} and continued for example in \cite{ambrosio}, and with estimation of discretely observed diffusion \cite{bibby}, where eigenfunctions play an important role. We will explore such links in future work. 
As a side note, we included in this paper a result showing an algebraic algorithm to calculate the canonical parameters $\theta$ from the expectation parameters $\eta$, which may be needed for implementing Equation \eqref{eq:etaode} in cases where one does not use eigenfunctions for $c$ and resorts to monomial sufficient statistics. It may be worth trying to find closed form solutions or approximations to move from expectation to canonical parameters for more general choices of the sufficient statistics $c$.  

Further important generalizations one may study are the  local and global projections of  Stochastic PDEs that extend the FPK PDE, for example in the filtering problem. In this case one has to manage carefully the stochastic part of the PDE. The SPDE is sensitive to the choice of stochastic calculus and the presence of noise, which cannot be differentiated, may compromise the local optimality of the projection, which cannot be assumed to hold pathwise a priori. In \cite{brigo98}, \cite{brigo99} and \cite{armstrongbrigomcss} a Stratonovich version of the infinite dimensional SPDE is considered and projected, but without any rigorous proof of local optimality. One then needs notions of optimality for SPDE that take into account the rough nature of the noise. This has been addressed for SDEs in \cite{armstrongBrigoproj} where optimal ways to project potentially high dimensional SDEs on low dimensional submanifolds have been derived, leading to two new projection methods, the Ito-vector and Ito-jet projections. Such methods have been formally extrapolated and briefly applied to the filtering SPDE in \cite{armstrongBrigoproj} and \cite{armstrongbrigoicms}, and are based on the jet bundle interpretation of Ito SDE's on manifolds \cite{armstrongBrigoJets}. In future work we will study rigorously the lift of the Ito vector and Ito jet projection to the infinite dimensional case of the filtering problem, based on replacing the high dimensional SDE to be projected with the infinite dimensional SPDE. 

Finally, one further potential application to filtering with continuous time observations is the following. Find particular systems for which the observation functions are among a set of eigenfunction of the state equation operator and where the optimal filter SPDE is approximated on a finite dimensional family based on said set of eigenfunctions. In this case, with exponential families and Hellinger distance the projection of the rough observations-driven part of the SPDE can be exact, see \cite{brigo99} for the case of exponential families. We could try to extend this result to simple mixtures and also combine it with the optimality in the projection of the drift part of the SPDE coming from generalizing the eigenfunctions result for the FPK equation, to see if we can obtain particularly good finite dimensional filters in this setting.  We might also explore the relationship between the discrete time observation case and the continuous time observations case, to see if the limit of the former can lead to the latter. This could be interesting also in cases where the state equation parameters are to be estimated and are not known \cite{schuppen}, as we find often in social sciences applications.


\begin{thebibliography}{99}

\bibitem{aggrawal74a} J. Aggrawal: Sur l'information de Fisher. In: Th\'eories de l'Information (J. Kampe de Feriet, ed.), Springer-Verlag, Berlin--New York 1974, pp. 111-117.

\bibitem{amari85a} Amari, S. Differential-geometrical methods in statistics, Lecture notes in statistics, Springer-Verlag, Berlin, 1985

%\bibitem{armstrongbrigoicms}
%John Armstrong and Damiano Brigo.
%\newblock Extrinsic projection of {I}t{\^o} {SDE}s on submanifolds with
%  applications to non-linear filtering.
%\newblock {\em To appear in: Nielsen, F., Critchley, F., \& Dodson, K. (Eds),
%  Computational Information Geometry for Image and Signal Processing, Springer
%  Verlag}, 2016.

\bibitem{ahmedbook}  Ahmed N. U. (1998). Linear and Nonlinear Filtering for Scientists and Engineers. World Scientific, Singapore. 

\bibitem{ahmed97} Ahmed, N. U. and S. Radaideh (1997). A powerful numerical technique solving the Zakai equation
for nonlinear filtering. Dynamics and Control 7 (3), 293--308.



\bibitem{ambrosio} L. Ambrosio, G. Savar\'e and L. Zambotti, Existence and stability for Fokker--Planck equations with log-concave reference measure, Probab. Theory Relat. Fields (2009) 145--517.



\bibitem{armstrongbrigomcss}
J.~Armstrong and D.~Brigo.
\newblock {N}onlinear filtering via stochastic {PDE} projection on mixture
  manifolds in {$L^2$} direct metric.
\newblock {\em Mathematics of Control, Signals and Systems}, 28(1):1--33, 2016.

\bibitem{armstrongBrigoJets}
J.~Armstrong and D.~Brigo.
\newblock Coordinate free stochastic differential equations as jets.
\newblock {\em http://arxiv.org/abs/1602.03931}, 2016.


\bibitem{armstrongBrigoproj}
J.~Armstrong and D.~Brigo.
\newblock Optimal approximation of SDEs on submanifolds: the Ito-vector and Ito-jet projections.
\newblock {\em http://arxiv.org/abs/1610.03887}, 2016.


\bibitem{armstrongbrigoicms}
J.~Armstrong and D.~Brigo.
\newblock Extrinsic projection of {I}t{\^o} {SDE}s on submanifolds with
  applications to non-linear filtering.
\newblock {\em In: Nielsen, F., Critchley, F., \& Dodson, K. (Eds),
  Computational Information Geometry for Image and Signal Processing, Springer
  Verlag}, 2016.



%\bibitem{crisan10} Bain, A., and Crisan, D. (2010). Fundamentals of Stochastic Filtering. Springer-Verlag, Heidelberg.

\bibitem{barndorff-nielsen78a} Barndorff-Nielsen, O.E. (1978). Information and Exponential Families. John Wiley and Sons, New York.


\bibitem{bibby} B. M. Bibby and M. Soerensen. Martingale estimation functions for discretely observed diffusion processes, Bernoulli, 1995, Vol. 1, N. 1-2, 17--39.

\bibitem{BorweinHuang95} Borwein, J.M., and Huang, W.Z. (1995).
A fast heuristic method for polynomial moment problems with
Boltzmann-Shannon entropy. {\em SIAM J. Optimization} 5, 68-99.


\bibitem{brigogyor}
Brigo, D.: On nonlinear {SDEs} whose densities evolve in a finite--dimensional
  family. In: Stochastic Differential and Difference Equations, Progress in
  Systems and Control Theory, vol.~23, pp. 11--19. Birkh{\"{a}}user Boston
  (1997)

\bibitem{brigoime}
{Brigo}, D., and Hanzon, B. {On some filtering problems arising in mathematical finance}.
  Insurance: {M}athematics and {E}conomics  22(1),  53--64 (1998)

%\bibitem{brigo99b} Brigo, D. Diffusion Processes, Manifolds of Exponential Densities, and Nonlinear Filtering, In: Ole E. Barndorff-Nielsen and Eva B. Vedel Jensen, editor, Geometry in Present Day Science, World Scientific, 1999

%\bibitem{brigo00}    Brigo, D, On SDEs with marginal laws evolving in finite-dimensional exponential families, STAT PROBABIL LETT, 2000, Vol: 49, Pages: 127 -- 134



\bibitem{brigo98}    Brigo, D, Hanzon, B, LeGland, F, A differential geometric approach to nonlinear filtering: The projection filter, IEEE T AUTOMAT CONTR, 1998, Vol: 43, Pages: 247 -- 252

\bibitem{brigo99}    Brigo, D, Hanzon, B, Le Gland, F, Approximate nonlinear filtering by projection on exponential manifolds of densities, BERNOULLI, 1999, Vol: 5, Pages: 495 -- 534

%\bibitem{Brigo2} D. Brigo, On the nice behaviour of the Gaussian
%Projection Filter with small observation noise,
%{\em Systems and Control Letters} {\bf 26} (1995)  363--370

%\bibitem{Brigo3} D. Brigo, New results on the Gaussian projection
%filter with small observation noise,  to appear in
%{\em Systems and Control Letters}.
%
\bibitem{BrigoPhD} D. Brigo, {\em Filtering by Projection on the Manifold
of Exponential Densities}, PhD Thesis, Free University of Amsterdam, 1996.

%\bibitem{BrigoCNR} D. Brigo,  On  diffusions
%with prescribed diffusion coefficients whose
%densities evolve in  prescribed exponential families,
%{\em Internal Report CNR--LADSEB} {\bf 02/96},
% CNR--LADSEB, Italy, March 1996.

%\bibitem{BrHaLe} D. Brigo, B. Hanzon, F. Le Gland,
% A differential geometric approach to nonlinear filtering:
%the projection filter, In {\em Proceedings of the Conference on Decision
%and Control}, IEEE-CSS, pp 4006--4011, New Orleans, 1995
%(extended version available on the internet at
%URL:~ftp://ftp.irisa.fr/techreports/1995/PI-914.ps.Z).

%\bibitem{ChalMaur} M. Chaleyat-Maurel, D. Michel,
%Des resultats de non-existence de filtre de dimension finie,
%{\em Stochastics}, {\bf 13} (1984) 83--102.

%\bibitem{daum86a} F. E. Daum, Exact Finite--Dimensional Nonlinear Filters,
%{\em IEEE Transactions on Automatic Control}, {\bf 31} (1986) 616--622.

\bibitem{brigopistone} Brigo, D., and Pistone, G. (2016).  Dimensionality reduction for measure-valued evolution equations in statistical manifolds.  In: Nielsen, F., Critchley, F., \& Dodson, K. (Eds),
{\em Proceedings of the conference on  Computational Information Geometry for Image and Signal Processing, Springer
  Verlag}, 2016. Earlier preprint version available at {\tt{http://arxiv.org/abs/1601.04189}}

\bibitem{brigopistone2} Brigo, D., and Pistone, G. (2016). Maximum likelihood eigenfunctions of the Fokker--Planck equation and Hellinger projection. Available at 
{\tt{http://arxiv.org/abs/1603.04348}}


%\bibitem{crisan} Crisan, D., and Rozovskii, B. (Eds) (2011). The Oxford Handbook of Nonlinear Filtering, Oxford University Press.

\bibitem{flandoli} G. Da Prato, F. Flandoli, M. R\"{o}ckner, Fokker-Planck Equations for SPDE with Non-trace-class Noise, Communications in Mathematics and Statistics, Volume 1, Issue 3,  281--304 (2013).


%\bibitem{davis81b} M. H. A. Davis, S. I. Marcus,
%An introduction to nonlinear filtering, in: M. Hazewinkel, J.
%C. Willems, Eds., {\em Stochastic Systems: The
%Mathematics of Filtering and Identification and Applications} (Reidel,
%Dordrecht, 1981) 53--75.

%\bibitem{FrieA} A. Friedman, {\em Stochastic Differential Equations and
%Applications},
%vol. I (Academic Press, New York, 1975).

\bibitem{girolami} M. Girolami and B. Calderhead, Riemann manifold Langevin and Hamiltonian Monte Carlo methods. Journal of the Royal Statistical Society: Series B (Statistical Methodology), 2011, 73, 123--214


%\bibitem{gyongy90a} I. Gy\"{o}ngy, N.V. Krylov, Stochastic partial differential
%equations with unbounded coefficients and applications II,
%{\em Stochastics}, {\bf 32} (1990) 165--180.


%\bibitem{elworthy82a} Elworthy, D. (1982). Stochastic Differential Equations on Manifolds. LMS Lecture Notes. 

\bibitem{hanzon87} Hanzon, B. A differential-geometric approach to approximate nonlinear filtering. In C.T.J. Dodson, Geometrization of Statistical Theory, pages 219 -- 223,ULMD Publications, University of Lancaster, 1987.

\bibitem{hut} Hanzon, B. and Hut, R. (1991). New results on the projection filter,  Proceedings of the First European Control Conference, Grenoble, 1991, Vol. I, pp. 623--628.

%\bibitem{hanzon89a} B. Hanzon, Identifiability, recursive identification and spaces of linear dynamical systems, CWI Tracts 63 and 64, CWI, Amsterdam, 1989

%\bibitem{HaMaSu} M. Hazewinkel, S.I.Marcus, and H.J. Sussmann,
%Nonexistence of finite dimensional filters for conditional statistics
%of the cubic sensor problem, {\em Systems and Control Letters} {\bf 3}
% (1983) 331--340.

%\bibitem{jacod87a} J. Jacod, A. N. Shiryaev, Limit theorems for stochastic processes. Grundlehren der Mathematischen Wissenschaften, vol. 288 (1987), Springer-Verlag, Berlin,

\bibitem{jazwinski70a} A. H. Jazwinski, {\em Stochastic Processes and
Filtering Theory}, Academic Press, New York, 1970.

\bibitem{otto} R. Jordan, D. Kinderlehrer and F. Otto,
The Variational Formulation of the Fokker--Planck Equation,
SIAM Journal on Mathematical Analysis, 1998, 
Vol. 29, N. 1, 1--17. 

%\bibitem{fujisaki72a} M. Fujisaki, G. Kallianpur, and H. Kunita (1972). Stochastic differential equations for the non linear filtering problem. 
%Osaka J. Math. Volume 9, Number 1 (1972), 19-40. 

\bibitem{kagan} Kagan, A.M. , Linnik, Y.V., and Rao, C.R. (1973). Characterization problems in Mathematical Statistics. John Wiley and Sons, New York.


%\bibitem{khasminskii} R. Z. Khasminskii (1980). Stochastic Stability of Differential Equations. Alphen aan den Reijn


\bibitem{kushner} H. Kushner, Approximations to optimal nonlinear filters. IEEE Trans. Automatic Control, 12, 546--556, 1967. 

%\bibitem{Levi} J. L\'evine, Finite dimensional realizations of
%stochastic PDE's and application to filtering,
%{\em Stochastics and Stochastic Reports} {\bf 43} (1991) 75--103.



\bibitem{liberzon} D. Liberzon and R. W. Brockett,
Spectral Analysis of Fokker--Planck and Related Operators Arising From Linear Stochastic Differential Equations,
SIAM Journal on Control and Optimization, 2000, 
Vol. 38, N. 5, 1453--1467


%\bibitem{LipShi} R.S. Liptser, A.N. Shiryayev, {\em Statistics of Random
%Processes I, General Theory} (Springer Verlag, Berlin, 1978).

\bibitem{maybeck} Maybeck, P. S. (1982). Stochastic models, estimation and
control. Academic Press.

%\bibitem{murray93a} M. Murray and J. Rice - Differential geometry and statistics, Monographs on Statistics and Applied Probability 48, Chapman and Hall, 1993.


\bibitem{naudts}
Naudts, J. (2008). Generalised Exponential Families and Associated Entropy
Functions,  Entropy, 10, 131-149.

\bibitem{Newton1}
N. J. Newton, An infinite-dimensional statistical manifold modelled on
  {H}ilbert space. J. Funct. Anal.  263(6),  1661--1681 (2012).
 
\bibitem{Newton2}
N.J. Newton, Infinite-dimensional manifolds of finite-entropy probability
  measures. In: Geometric science of information, Lecture Notes in Comput.
  Sci., vol. 8085, pp. 713--720. Springer, Heidelberg (2013)

\bibitem{Newton3}
N. J. Newton, Information geometric nonlinear filtering. Infin. Dimens. Anal.
  Quantum Probab. Relat. Top.  18(2),  1550014, 24 (2015).


%\bibitem{kenney} Kenney, J., Stirling, W. Nonlinear Filtering of Convex Sets of Probability Distributions. Presented at the 1st International Symposium on Imprecise Probabilities and Their Applications, Ghent, Belgium, 29 June - 2 July 1999

%\bibitem{OconPard} D. Ocone, E. Pardoux,
%{\em A Lie algebraic criterion for non-existence of finite dimensionally
%computable filters}, Lecture notes in mathematics 1390, 197--204
%(Springer Verlag, 1989)




\bibitem{pavliotis}
G.~A. Pavliotis.
\newblock {\em Stochastic Processes and Applications: Diffusion Processes, the
  Fokker-Planck and Langevin Equations}.
\newblock Springer, Heidelberg, 2014.


\bibitem{pistonesempi} Pistone, G., and Sempi, C. (1995). An Infinite Dimensional Geometric Structure On the
space of All the Probability Measures Equivalent to a Given one. The Annals of Statistics
23(5), 1995


\bibitem{santacroce}
Santacroce, M., Siri, P., Trivellato, B., New results on mixture and
  exponential models by {O}rlicz spaces, Bernoulli,  
  Volume 22, Number 3 (2016), 1431-1447.

\bibitem{schuppen} van Schuppen, J. H. (1983). Convergence results for continuous-time adaptive stochastic filtering algorithms, J. Math. Anal. Appl. 96, 209--225.


\end{thebibliography}
\end{document}